\documentstyle[amscd,amssymb,verbatim,epsf]{amsart}

\begin{document}
\theoremstyle{plain}
\newtheorem{Thm}{Theorem}
\newtheorem{Cor}{Corollary}
\newtheorem{Con}{Conjecture}
\newtheorem{Main}{Main Theorem}
\newtheorem{Lem}{Lemma}
\newtheorem{Prop}{Proposition}

\theoremstyle{definition}
\newtheorem{Def}{Definition}
\newtheorem{Note}{Note}
\newtheorem{Ex}{Example}

\theoremstyle{remark}
\newtheorem{notation}{Notation}
\renewcommand{\thenotation}{}

\errorcontextlines=0
\numberwithin{equation}{section}
\renewcommand{\rm}{\normalshape}%

\title
   {Reflection in a Translation Invariant Surface}

\author{Brendan Guilfoyle}
\address{Brendan Guilfoyle\\
          Department of Mathematics and Computing \\
          Institute of Technology, Tralee \\
          Clash \\
          Tralee  \\
          Co. Kerry \\
          Ireland.}
\email{brendan.guilfoyle@@ittralee.ie}

\author{Wilhelm Klingenberg}
\address{Wilhelm Klingenberg\\
 Department of Mathematical Sciences\\
 University of Durham\\
 Durham DH1 3LE\\
 United Kingdom.}
\email{wilhelm.klingenberg@@durham.ac.uk }

\date{September 26, 2005}

\begin{abstract}
We prove that the focal set generated by the reflection of a point
source off a translation invariant surface consists of two sets: a
curve and a surface. The focal curve lies in the plane orthogonal to
the symmetry direction containing the source, while the focal surface is
translation invariant. 

This is done by constructing explicitly the focal set of the
reflected line congruence (2-parameter family of oriented lines in
${\Bbb{R}}^3$) with the aid of the natural complex structure on the space
of all oriented affine lines.
\end{abstract}

\maketitle

The purpose of this paper is to prove the following Theorem:
\vspace{0.1in}

\noindent{\bf Main Theorem}:

{\it The focal set generated by the reflection of a point source off a
translation invariant surface consists of two sets: a curve and a
surface. The focal curve lies in the plane orthogonal to the symmetry
direction containing the source, while the focal surface is
translation invariant. 
}
\vspace{0.1in}

In contrast to the focal surface, the reflected wavefront is
not translation invariant, in general. 

There have been many investigations of generic focal sets of line
congruences \cite{agv} \cite{bgg} \cite{ist}. Rather than work in the
generic setting, we compute the focal set explicitly in this special case.
This we do by applying recent work on immersed
surfaces in the space ${\Bbb{T}}$  of oriented affine lines in
${\Bbb{R}}^3$ \cite{gak2} \cite{gak3}.  The next section contains a
summary of the background material on the complex geometry of ${\Bbb{T}}$
and the focal sets of arbitrary line congruences. It also details the
reflection of a line congruence in an oriented surface in ${\Bbb{R}}^3$.

In Section 2 we solve the problem of reflection
of a point source off an arbitrary translation invariant surface
(Proposition \ref{p:transymref}). We then compute the focal set and
thus prove the Main Theorem.

\section{The Parametric Approach to Geometric Optics}

Let ${\Bbb{T}}$ be the set of oriented affine lines in Euclidean
${\Bbb{R}}^3$, which, by parallel translation, can be identified with the
tangent bundle to the 2-sphere, $\mbox{T}{\Bbb{P}}^1$. We summarise
now the main features of this identification (further 
details can be found in \cite{gak2} \cite{gak3}).

\newpage

The canonical projection $\pi:{\Bbb{T}}\rightarrow {\Bbb{P}}^1$
assigns an oriented line, or {\it ray},  to its direction, and we also
have the double fibration:  

\vspace{0.1in}
\begin{center}

\unitlength0.5cm

\begin{picture}(12,9)
\put(-2.6,7.5){${\Bbb{T}}\times{\Bbb{R}}$}
\put(-1,7){\vector(1,-1){2.8}}
\put(0.5,5.8){$\Phi$}
\put(-2.2,3){${\Bbb{T}}$}
\put(-2,7){\vector(0,-1){3}}
\put(2,3.4){${\Bbb{R}}^3$}
\put(7,0){\mbox{\epsfbox{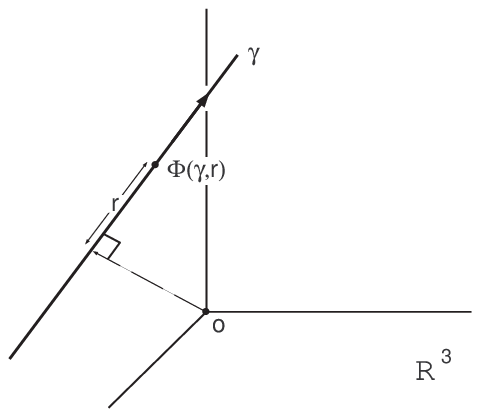}}}
\end{picture}

\end{center}

\vspace{0.1in}

\noindent In the diagram, the lefthand map is projection onto the
first factor. The 
mapping on the right, denoted by $\Phi$, takes 
$(\gamma,r)\in{\Bbb{T}}\times{\Bbb{R}}$ to the point on the oriented line
$\gamma$ in ${\Bbb{R}}^3$ that lies an affine parameter distance $r$
from the point on $\gamma$ closest to the origin (as shown).

Let $\xi$ be the local coordinates on ${\Bbb{P}}^1$ obtained by
stereographic projection from the south pole. This can be extended to
coordinates $(\xi,\eta)$ on ${\Bbb{T}}$ minus the fibre over the south
pole. The map $(\xi,\eta,r)\mapsto\Phi(\xi,\eta,r)=(z(\xi,\eta,r),t(\xi,\eta,r))$ has the following cooordinate expression \cite{gak2}:
\begin{equation}\label{e:coord}
z=\frac{2(\eta-\overline{\eta}\xi^2)+2\xi(1+\xi\overline{\xi})r}{(1+\xi\overline{\xi})^2}
\qquad\qquad
t=\frac{-2(\eta\overline{\xi}+\overline{\eta}\xi)+(1-\xi^2\overline{\xi}^2)r}{(1+\xi\overline{\xi})^2},
\end{equation}
where $z=x^1+ix^2$, $t=x^3$ and ($x^1$, $x^2$, $x^3$) are Euclidean
coordinates in ${\Bbb{R}}^3$.

\begin{Def}
A {\it line congruence} is an immersed surface
$f:\Sigma\rightarrow{\Bbb{T}}$, i.e. a 2-parameter family of
oriented lines in ${\Bbb{R}}^3$. A {\it smoothly parameterised} line
congruence is a smoothly immersed surface $f:\Sigma\rightarrow{\Bbb{T}}$
together with an open cover $\{U_\alpha\}$ of $\Sigma$ and diffeomorphisms
${\Bbb{C}}\rightarrow U_\alpha:\mu_\alpha\mapsto\gamma$. For short we 
denote a parameterisation simply by $\mu$, and assume that all maps
are at least $C^1$-smooth.

\end{Def}

The first order
properties of such a family can described by two complex functions,
the {\it optical scalars}: 
$\rho,\sigma:\Sigma\times{\Bbb{R}}\rightarrow{\Bbb{C}}$, which are
defined relative to an orthonormal frame in ${\Bbb{R}}^3$ adapted to
the congruence. The real part $\theta$ and the imaginary part
$\lambda$ of $\rho$ are the {\it divergence} and {\it twist} of the
congruence, while $\sigma$ is the {\it shear} \cite{par}. 

For a parameterised line congruence we compose with the coordinates
above to get $\mu\mapsto(\xi(\mu,\bar{\mu}),\eta(\mu,\bar{\mu}))$. The
optical scalars then, with a natural choice of orthonormal frame, have the
following expressions \cite{gak2}:  
\begin{equation}\label{e:spinco}
\rho=\theta+\lambda i=\frac{ \partial^+\eta\overline{\partial}\;\overline{\xi} -\partial^-\eta\partial\overline{\xi}}
{\partial^-\eta\overline{\partial^-\eta}-\partial^+\eta\overline{\partial^+\eta}}
\qquad\qquad
\sigma=\frac{\overline{\partial^+\eta}\partial\overline{\xi} -\overline{\partial^-\eta}\;\overline{\partial}\;\overline{\xi}}
{\partial^-\eta\overline{\partial^-\eta}-\partial^+\eta\overline{\partial^+\eta}},
\end{equation}
where
\[
\partial^+\eta\equiv\partial \eta+r\partial\xi-\frac{2\eta\overline{\xi}\partial \xi}{1+\xi\overline{\xi}}
\qquad\qquad
\partial^-\eta\equiv\overline{\partial} \eta+r\overline{\partial}\xi-\frac{2\eta\overline{\xi}\;\overline{\partial} \xi}{1+\xi\overline{\xi}},
\]
and $\partial$ and $\bar{\partial}$ are differentiation with respect
to $\mu$ and $\bar{\mu}$, respectively.

\begin{Def}
The {\it curvature} of a line congruence is defined to be
$\kappa=\rho\bar{\rho}-\sigma\bar{\sigma}$. A line congruence is {\it flat}
if $\kappa=0$. A line congruence ${\Sigma}\subset{\Bbb{T}}$ is flat
iff the rank of the projection $\pi:{\Bbb{T}}\rightarrow {\Bbb{P}}^1$
restricted to $\Sigma$ is non-maximal.
\end{Def}

\begin{Def}
A point $p$ on a line $\gamma$ in a line congruence is a
{\it focal point} if $\rho$ and $\sigma$ blow-up at $p$. The set of
focal points of a line congruence $\Sigma$ generically form surfaces
in ${\Bbb{R}}^3$, which will be referred to as the {\it focal
surfaces} of $\Sigma$. 
\end{Def}

\begin{Thm}\label{t:focs}

The focal set of a parametric line congruence $\Sigma$ is
\[
\{\Phi(\gamma,r)|\;\gamma\in\Sigma\;\;\mbox{and}\;\;
1-2\theta_0r+(\rho_0\overline{\rho}_0-\sigma_0\overline{\sigma}_0)r^2=0\}, 
\]
where the coefficients of the quadratic equation are given by (\ref{e:spinco})
at $r=0$. 
\end{Thm}
\begin{pf}
In terms of the affine
parameter $r$ along a given line, the Sachs equations, which $\sigma$ and
$\rho$ must satisfy, are \cite{par}:
\[
\frac{\partial \rho}{\partial r}=\rho^2+\sigma\overline{\sigma}
\qquad 
\frac{\partial \sigma}{\partial r}=(\rho+\overline{\rho})\sigma.
\]
These are equivalent to the vanishing of certain components of the
Ricci tensor of the Euclidean metric. They have solution:
\[
\rho=\frac{\rho_0-(\rho_0\overline{\rho}_0-\sigma_0\overline{\sigma}_0)r}
  {1-2\theta_0r+(\rho_0\overline{\rho}_0
                              -\sigma_0\overline{\sigma}_0)r^2}
\qquad
\sigma=\frac{\sigma_0}
  {1-2\theta_0r+(\rho_0\overline{\rho}_0-\sigma_0\overline{\sigma}_0)r^2},
\]
where  $\sigma_0$, $\theta_0$ and $\rho_0$ are the values of the optical
scalars at $r=0$. The theorem follows.
\end{pf}

This has the following corollary:

\begin{Cor}\label{c:quad}

Let $\Sigma$ be a line congruence, $\rho=\theta+\lambda i$,
$\sigma$ the associated optical scalars and $\rho_0$, $\theta_0$,
$\lambda_0$, $\sigma_0$ their values at $r=0$. 

If $\Sigma$ is flat with non-zero divergence, then there
exists a unique focal surface $S$ given by
$r=(2\theta_0)^{-1}$. If it is flat with zero divergence,
then the focal set is empty.

If $\Sigma$ is non-flat, then there exists a unique focal point on each
line iff  $|\sigma_0|^2=\lambda_0^2$, there exist two focal
points on each line iff $|\sigma_0|^2<\lambda_0^2$ and there are no
focal points on each line iff $|\sigma_0|^2>\lambda_0^2$. The
focal set is given by
\[
r=\frac{\theta_0\pm(|\sigma_0|^2 -\lambda_0^2)^{\frac{1}{2}}}
            {\rho_0\bar{\rho}_0-\sigma_0\bar{\sigma}_0}.
\]

\end{Cor}

\begin{pf}
The focal set of a parameterised line congruence
are given by $r=r(\mu,\bar{\mu})$ satisfying the quadratic equation
in Theorem \ref{t:focs}. If $\kappa=0$, then there is none or one
solution depending on whether $\theta_0=0$ or not.

If $\kappa\neq 0$ then there are two, one or no solutions iff
$|\sigma_0|^2-\lambda_0^2$ is greater than, equal to or less than
zero (respectively). 

The solution of the quadratic equation in each case is as stated.
\end{pf}

Given a line congruence $f:\Sigma\rightarrow{\Bbb{T}}$, a map
$s:\Sigma\rightarrow {\Bbb{R}}$ determines a map
$\Sigma\rightarrow{\Bbb{R}}^3$ by
$\gamma\mapsto\Phi(f(\gamma),s(\gamma))$ for $\gamma\in\Sigma$.
With a local parameterisation $\mu$ of $\Sigma$, we get a
map ${\Bbb{C}}\rightarrow{\Bbb{R}}^3$ which comes from 
substituting $r=s(\mu,\bar{\mu})$ in equations (\ref{e:coord}).

Of particular interest are the surfaces in ${\Bbb{R}}^3$ orthogonal to
the line congruence - when the line congruence is {\it normal}. These
exist iff the twist of the congruence vanishes, and the surfaces are
obtained from the solutions of the following equation \cite{gak2}:
\begin{equation}\label{e:intsur}
\bar{\partial} r=\frac{2\eta\bar{\partial}\bar{\xi}+2\bar{\eta}\bar{\partial}\xi}{(1+\xi\bar{\xi})^2}.
\end{equation}

We turn now to the reflection of an oriented line in a surface in
${\Bbb{R}}^3$. This is equivalent to the action of a certain
group on the space of oriented lines, as described by \cite{gak3}:

\begin{Thm}\label{t:ref}
Consider a parametric line congruence
$\xi=\xi_1(\mu_1,\bar{\mu}_1)$, $\eta=\eta_1(\mu_1,\bar{\mu}_1)$
reflected off an oriented surface with parameterised normal line
congruence $\xi=\xi_0(\mu_0,\bar{\mu}_0)$,
$\eta=\eta_0(\mu_0,\bar{\mu}_0)$ and $r=r_0(\mu_0,\bar{\mu}_0)$
satisfying (\ref{e:intsur}) with $\xi=\xi_0$ and $\eta=\eta_0$. Then
the reflected line congruence is 
\begin{equation}\label{e:reflaw1}
\xi=\frac{2\xi_0\bar{\xi}_1+1-\xi_0\bar{\xi}_0}
           {(1-\xi_0\bar{\xi}_0)\bar{\xi}_1-2\bar{\xi}_0},
\end{equation}
\begin{equation}\label{e:reflaw2}
\eta={\textstyle \frac{(\bar{\xi}_0-\bar{\xi}_1)^2}
         {((1-\xi_0\bar{\xi}_0)\bar{\xi}_1-2\bar{\xi}_0)^2}}\eta_0
       -{\textstyle\frac{(1+\xi_0\bar{\xi}_1)^2}
         {((1-\xi_0\bar{\xi}_0)\bar{\xi}_1-2\bar{\xi}_0)^2}}\bar{\eta}_0
+{\textstyle\frac{(\bar{\xi}_0-\bar{\xi}_1)(1+\xi_0\bar{\xi}_1)(1+\xi_0\bar{\xi}_0)}
         {((1-\xi_0\bar{\xi}_0)\bar{\xi}_1-2\bar{\xi}_0)^2}}r_0,
\end{equation}
where the incoming rays are only reflected if they satisfy the
intersection equation 
\begin{equation}\label{e:int}
\eta_1={\textstyle \frac{(1+\bar{\xi}_0\xi_1)^2}{(1+\xi_0\bar{\xi}_0)^2}}\eta_0
-{\textstyle \frac{(\xi_0-\xi_1)^2}{(1+\xi_0\bar{\xi}_0)^2}}\bar{\eta}_0
+{\textstyle \frac{(\xi_0-\xi_1)(1+\bar{\xi}_0\xi_1)}{1+\xi_0\bar{\xi}_0}}r_0.
\end{equation}
\end{Thm}

\vspace{0.2in}
By virtue of the intersection equation, an alternative way of writing
(\ref{e:reflaw2}) is 
\begin{equation}\label{e:reflaw3}
\eta={\textstyle \frac{-(1+\xi_0\bar{\xi}_0)^2}{((1-\xi_0\bar{\xi}_0)\bar{\xi}_1-2\bar{\xi}_0)^2}}\bar{\eta}_1
+{\textstyle \frac{2(\bar{\xi}_0-\bar{\xi}_1)(1+\xi_0\bar{\xi}_1)(1+\xi_0\bar{\xi}_0)}{((1-\xi_0\bar{\xi}_0)\bar{\xi}_1-2\bar{\xi}_0)^2}}r_0.
\end{equation}
The geometric content of this is: reflection of an oriented line
can be decomposed into a sum of rotation about the origin (the
derived action of PSL(2,${\Bbb{C}}$) on $\mbox{T}{\Bbb{P}}^1$) and
translation (a fibre-mapping on $\mbox{T}{\Bbb{P}}^1$).

\section{Reflection off a translation invariant surface}

Consider a translation invariant surface with axis
lying along the $x^3-$axis in ${\Bbb{R}}^3$. Such a cylinder can be
parametrised by $(u,v)\mapsto(z=z_0(u),t=v)$ for
$(u,v)\in{\Bbb{R}}^2$, where $z=x^1+ix^2$, $t=x^3$ and ($x^1$, $x^2$,
$x^3$) are Euclidean coordinates in ${\Bbb{R}}^3$. 

\begin{Prop}
The normal congruence of a translation invariant surface is:
\begin{equation}\label{e:surdir}
\xi_0=\pm\left[-\frac{\dot{z}_0}{\dot{\bar{z}}_0}\right]^{\frac{1}{2}}
\end{equation}
\begin{equation}\label{e:surperp}
\eta_0={\textstyle \frac{1}{2}}(z_0-2v\xi_0-\bar{z}_0\xi_0^2)
\qquad\qquad
r_0=\frac{\bar{\xi}_0z_0+\xi_0\bar{z}_0+(1-\xi_0\bar{\xi}_0)v}{1+\xi_0\bar{\xi}_0},
\end{equation}
where a dot represents differentiation with respect to $u$ and the
choice of sign of $\xi_0$ is one of orientation.
\end{Prop}
\begin{pf}
Let $\xi_0\in S^2$ be the direction of the normal. This corresponds to
the vector
\[
\frac{2\xi_0}{1+\xi_0\bar{\xi}_0}\frac{\partial}{\partial z} +\frac{2\bar{\xi}_0}{1+\xi_0\bar{\xi}_0}\frac{\partial}{\partial \bar{z}}+\frac{1-\xi_0\bar{\xi}_0}{1+\xi_0\bar{\xi}_0}\frac{\partial}{\partial t}
\]
The vanishing of the inner product of this vector with the push
forward of $\frac{\partial}{\partial u}$ and $\frac{\partial}{\partial
  v}$ yields
\[
1-\xi_0\bar{\xi}_0=0 \qquad\qquad\qquad \dot{z}_0\bar{\xi}_0
    +\dot{\bar{z}}_0\xi_0=0
\]
and equations (\ref{e:surdir}) follows.

Equations (\ref{e:surperp}) come from inverting
equation (\ref{e:coord}) for $\eta$ and $r$.
\end{pf}

We consider the reflection of a point source off a surface that is
translation invariant along the $x^3$-axis. By a translation of
surface and source we move the point source to the origin. The line
congruence consisting of all oriented lines through the origin is
given by $\eta_1=0$ and $\xi_1\in S^2$.

\begin{Prop}\label{p:transymref}
The reflection of a point source at the origin off a translation
invariant surface is $\xi=-\xi_0^2\bar{\xi}_1$ and
$\eta=\left(\bar{\xi}_0-\xi_0\bar{\xi}_1^2\right)\xi_0^2r_0$ where
\begin{equation}\label{e:inteq}
\xi_1=\frac{u\pm(v^2+z_0\bar{z}_0)^{\frac{1}{2}}}{\bar{z}_0}
\end{equation}
and $\xi_0$ is given by (\ref{e:surdir}). Here the $\pm$ refers to the
two oriented lines from the source that intersect any given point in
${\Bbb{R}}^3$, and is chosen so that the ray goes from the source to
the point of reflection.
\end{Prop}
\begin{pf}
The reflection equations contained in Theorem
\ref{t:ref}, with $\eta_1=0$, $\xi_1\in S^2$  yield the stated
reflected line congruence and the intersection equation has solution
(\ref{e:inteq}). 
\end{pf}

\vspace{0.1in}
\noindent{\bf Main Theorem}:

{\it The focal set generated by the reflection of a point source off a
translation invariant surface consists of two sets: a curve and a
surface. The focal curve lies in the plane orthogonal to the symmetry
direction containing the source, while the focal surface is
translation invariant. 
}
\vspace{0.1in}

\begin{pf}
The local parameter we choose is $\mu=u+iv$ and the 
spin coefficients are computed by inserting the reflected line
congruence in Proposition \ref{p:transymref} into equation
(\ref{e:spinco}). The focal set is then determined by inserting the
reflected line congruence and the solutions $r$ of the quadratic
equation in Theorem \ref{t:focs} into (\ref{e:coord}). Assuming that
both $\pm$ in (\ref{e:surdir}) and (\ref{e:inteq}) are positive, we
obtain the following focal set: 
\[
z=\frac{z_0\dot{\bar{z}}_0-\bar{z}_0\dot{z}_0}{\dot{\bar{z}}_0}
\qquad\qquad\qquad
t=0,
\]
and 
\[
z=\frac{2\ddot{z}_0\dot{\bar{z}}_0z_0^2\bar{z}_0
        -2\ddot{\bar{z}}_0\dot{z}_0z_0^2\bar{z}_0
        +\dot{z}_0^3\bar{z}_0^2-2\dot{z}_0^2\dot{\bar{z}}_0z_0\bar{z}_0
        +\dot{z}_0\dot{\bar{z}}_0^2z_0^2}
       {2\ddot{z}_0\dot{\bar{z}}_0z_0\bar{z}_0
       -2\ddot{\bar{z}}_0\dot{z}_0z_0\bar{z}_0
       -\dot{z}_0^2\dot{\bar{z}}_0\bar{z}_0
       +\dot{z}_0\dot{\bar{z}}_0^2z_0},
\]
\[
t=\frac{2v\;z_0\bar{z}_0(\ddot{\bar{z}}_0\dot{z}_0-\ddot{z}_0\dot{\bar{z}}_0)}
       {2\ddot{z}_0\dot{\bar{z}}_0z_0\bar{z}_0
       -2\ddot{\bar{z}}_0\dot{z}_0z_0\bar{z}_0
       -\dot{z}_0^2\dot{\bar{z}}_0\bar{z}_0
       +\dot{z}_0\dot{\bar{z}}_0^2z_0}.
\]
The first of these is a curve in the $x^1x^2$-plane parameterised by
$u$, and the second is a surface that is invariant in the
$x^3$-direction.

A similar result holds for the other choices of orientation.
\end{pf}


\begin{thebibliography}{10}
\bibitem{agv}
 V. I. Arnold, S. M. Gusein-Zade, A. N. Varchenko, Singularities of
differentiable maps, Volume 1, Birkhaeuser, Basel, 1986.
\bibitem{bgg}
J. Bruce, P. Giblin and C. Gibson, {\it On caustics by reflection}, Topology {\bf 21} (1982), 179-199.
\bibitem{gak2}
B. Guilfoyle and W. Klingenberg, {\it Generalised surfaces in
  ${\Bbb{R}}^3$}, Math. Proc. R. Ir. Acad. {\bf 104A(2)}, 199-209.
\bibitem{gak3}
B. Guilfoyle and W. Klingenberg, {\it Reflection of a wave off a
  surface}, Journal of Geometry (to appear) [math.DG/0406212]
\bibitem{ist}
S. Izumiya, K. Saji and N. Takeuchi, {\it Singularities of line
  congruences}, Proceedings of the Royal Society of Edinburgh {\bf
  133A} (2003), 1341-1359. 
\bibitem{par}
R. Penrose and W. Rindler, Spinors and spacetime, Volume 1 and 2, Cambridge
University Press, Cambridge, 1986.








\end{thebibliography}
\end{document}